\documentclass{amsart}
\usepackage{amsfonts,amscd,amssymb,amsmath,mathrsfs}
\newtheorem{theorem}{Theorem}[section]
\newtheorem{definition}{Definition}[section]
\newtheorem{lemma}[theorem]{Lemma}
\newtheorem{proposition}[theorem]{Proposition}
\newtheorem{corollary}[theorem]{Corollary}
\theoremstyle{remark}
\newtheorem{rmk}{Remark}[section]
\theoremstyle{remark}
\newtheorem{example}{Example}[section]
\DeclareMathOperator{\rank}{rank}

\DeclareMathOperator{\End}{End}
\numberwithin{equation}{section}
\newcommand{\module}{\operatorname{mod}}
\newcommand{\La}{\operatorname{\Lambda}}
\newcommand\D{\operatorname{D}}
\newcommand\Img{\operatorname{Im}}
\newcommand{\Hom}{\operatorname{Hom}}
\newcommand{\Ker}{\operatorname{Ker}}
\newcommand\Coker{\operatorname{Coker}}

\begin{document}
\title[Preprojective Algebras and 
Multiplication Maps of Maximal Rank]{Almost Split Morphisms, Preprojective Algebras and 
Multiplication Maps of Maximal Rank}
\author{Steven P. Diaz}
\address{Department of Mathematics, Syracuse University, Syracuse,
 New York 13244-1150, United States of America}
\email{spdiaz@syr.edu}
\author{Mark Kleiner}
\email{mkleiner@syr.edu}
\thanks{The second-named author is partially supported by a grant from NSA}
\maketitle

\begin{abstract}
With a grading previously 
introduced by the second-named author, the multiplication maps in the preprojective algebra 
satisfy a maximal rank property that is similar to the maximal rank 
property proven by Hochster and Laksov for the multiplication maps in 
the commutative polynomial ring. The result follows from a more 
general theorem about the maximal rank property of a  minimal almost split 
morphism, which also yields a quadratic inequality for the dimensions of indecomposable modules involved.
\end{abstract}

\section{Introduction}
Let $k$ be a  field and let $R$ be the polynomial ring in $n$ commuting variables over $k$. Let $R_i$ be its $i^{th}$ graded piece consisting of 
homogeneous polynomials of degree $i$. A result of Hochster and 
Laksov \cite{HL1} says that if $i\ge2$ and $\mathnormal{V}\subset R_i$ is a 
general subspace then the natural multiplication map from 
$\mathnormal{V}\otimes R_1$ to $R_{i+1}$ has maximal rank, that is is 
either injective or surjective, and it is not known what happens if one replaces $R_1$ by $R_d$ for $d>1$. One may wonder which other graded 
rings have a similar property. 

In \cite{K1} a 
new grading on the preprojective algebra was introduced. In this paper we show that with this grading,
the preprojective algebra of a finite quiver without oriented cycles satisfies a property 
analogous to the Hochster-Laksov property for polynomial rings, and much 
of our proof is quite similar to their proof. At 
one point the proof for preprojective algebras becomes easier than 
the proof for polynomial rings: some of the more complicated 
dimension counts needed for polynomial rings are not needed for 
preprojective algebras. This allows us to obtain a result for 
preprojective algebras that is stronger than the analogous result for 
polynomial rings.

The key to making things work is the fact that the multiplication-by-arrow maps into a fixed homogeneous component of the infinite dimensional (in general) preprojective algebra  give rise to a minimal right almost split morphism of modules over the finite dimensional path algebra of the quiver \cite{K1}, which implies the maximal rank property. In fact we show that a  minimal right almost split morphism $g:B\to C$ of finite dimensional modules over a $k$-algebra satisfies a maximal rank property analogous to the Hochster-Laksov property for polynomial rings, and if $C$  is not projective and $B_1,\dots,B_l$ are the nonisomorphic indecomposable summands of $B$ then $\dim_k C<(\dim_k B_1)^2+\dots+(\dim_k B_l)^2$.  We do not know what happens if multiplication by arrows is replaced by multiplication by paths of fixed length greater than one.

There is a natural dual to the Hochster-Laksov maximal rank property, and the two properties always occur simultaneously.  We give two explanations of this fact, one general  homological and the other based on the vector space duality $\D=\Hom_k(\ ,k)$.  As a consequence, the multiplication-by-arrow maps out of a fixed homogeneous component of the preprojective algebra satisfy the dual Hochster-Laksov maximal rank property. 

The 
organization of the paper is as follows. In Section 2 we prove a 
theorem that gives a general situation in which one can obtain a 
maximal rank property analogous to the Hochster-Laksov property for 
polynomial rings. This general situation does not include the 
polynomial ring as a special case. In Section 3 we review some facts 
about almost split morphisms and preprojective algebras and then show 
that almost split morhisms in general and the preprojective algebra 
in particular fit into the general set up of Section 2. In Section 4 
we use the material in Sections 2 and 3 to obtain results for the 
preprojective algebra that look very analogous to the Hochster-Laksov 
result for polynomial rings. Then we conclude with some examples to 
illustrate the results.

In this paper for simplicity we work over a fixed algebraically closed field $k$ and $\dim$ always means $\dim_k$.  For unexplained terminology we refer the reader to ~\cite{ars}.

\section{The General Theorem}

Let $\mathnormal{V_1,V_2,...,V_l,W_1,W_2,...,W_l,U}$ be finite 
dimensional vector spaces. Let $\mathnormal{T}$ be a  
linear transformation from the direct sum of the tensor products 
$\mathnormal{V_1\otimes W_1}$, $\mathnormal{V_2\otimes W_2}$, ... ,$\mathnormal{V_l\otimes W_l}$ to 
$\mathnormal{U}$.
\begin{equation}\label{mapT}\mathnormal{T}:\bigoplus_{i=1}^{\mathnormal{l}}(\mathnormal{V_i\otimes 
W_i})\to\mathnormal{U}
\end{equation}
\begin{definition} 
\label{D1:uni}
We say that $\mathnormal{T}$ satisfies the right 
omnipresent maximal rank property if and only if for every choice of 
subspaces $\mathnormal{W'_i}\subset\mathnormal{W_i}$ for $i=1,...,l$ 
the restriction of $\mathnormal{T}$ to the direct sum of the tensor 
products $V_1\otimes W'_1,V_2\otimes W'_2,...$, $V_l\otimes 
W'_l$ has maximal rank, that is, is either injective or 
surjective.
\end{definition}

Notice that if ${T}$ satisfies the right omnipresent maximal rank property then so does its restriction to $\oplus^l_{i=1}V_i\otimes W'_i$, and ${T}$ must itself have maximal rank.  Every injective $\mathnormal{T}$ satisfies the right omnipresent maximal rank property. The interesting case is when $\mathnormal{T}$ is surjective but not injective.

Denote by $\End{(\mathnormal{V}_i)}$ 
the $k$-algebra of linear operators on $\mathnormal{V}_i$. The tensor product $\mathnormal{V}_i\otimes \mathnormal{W}_i$ is a left $\End{(\mathnormal{V}_i)}$-module by means of $\varphi_i\cdot (v_i\otimes 
w_i)=\varphi_i(v_i)\otimes w_i, \varphi_i\in\End{(\mathnormal{V}_i)},v_i\in V_i, w_i\in W_i$. Applying this to each term of the 
direct sum one obtains a bilinear evaluation 
map
\[e:\prod_{i=1}^l\End{(\mathnormal{V}_i)}\times\bigoplus_{i=1}^{\mathnormal{l}}(\mathnormal{V_i\otimes 
W_i})\to\bigoplus_{i=1}^{\mathnormal{l}}(\mathnormal{V_i\otimes 
W_i}).\]
Denote $\prod_{i=1}^l\End{(\mathnormal{V}_i)}$ by 
$\mathnormal{B}$. The map $e$ defines a structure of a left $B$-module on $\bigoplus_{i=1}^{\mathnormal{l}}(\mathnormal{V_i\otimes 
W_i})$. Notice that $\mathnormal{B}$ has dimension 
$\Sigma(\dim{\mathnormal{V_i}})^2$. Let $\mathnormal{P}_i$ be the 
projective space of one dimensional subspaces of 
$\mathnormal{V}_i\otimes\mathnormal{W}_i$ and let $\mathnormal{P}$ be 
the projective space of one dimensional subspaces of 
$\bigoplus_{i=1}^{\mathnormal{l}}(\mathnormal{V_i\otimes W_i})$. 
Notice that $\mathnormal{P}$ has dimension 
$\Sigma(\dim{\mathnormal{V_i}}\dim{\mathnormal{W_i}})-1$. We shall 
study the product $\mathnormal{B}\times\mathnormal{P}$ together with 
its two projection maps $\pi_1$ onto $\mathnormal{B}$ and $\pi_2$ 
onto $\mathnormal{P}$. 

Since the evaluation map $e$ is bilinear, we may conclude that the inverse image under $e$ of $\Ker T$ is a Zariski closed subset of the domain of $e$. Furthermore using bilinearity again we see that $e^{-1}(\Ker T)$ is the
affine cone over a Zariski closed subset of 
$\mathnormal{B}\times\mathnormal{P}$. We denote this subset by 
$\mathnormal{Y}$. For each $i$ from $1$ to $l$ let $\mathnormal{X_i}$ 
be an irreducible quasiprojective subset of $\mathnormal{P_i}$ and 
let $\mathnormal{C(X_i)}$ be  its corresponding affine cone in 
$\mathnormal{V}_i\otimes\mathnormal{W}_i$. Let $\mathnormal{X}$ be 
the irreducible quasiprojective subset of $\mathnormal{P}$ 
corresponding to 
$\mathnormal{C(X_1)}\times\mathnormal{C(X_2)}\times...\times\mathnormal{C(X_l)}$. 
Notice that 
$\dim{\mathnormal{X}}=\sum_{i=1}^{l}\dim\mathnormal{C(X_i)}-1$.

\begin{theorem}\label{T:main}
Assume 
that $\mathnormal{T}$ satisfies the right omnipresent maximal rank 
property and that 
$\sum_{i=1}^{l}\dim\mathnormal{C(X_i)}\leq\dim\mathnormal{U}$. Then 
$\pi_1(\pi_2^{-1}(\mathnormal{X})\cap\mathnormal{Y})$ is contained in 
a proper Zariski closed subset of $\mathnormal{B}$.
\end{theorem}

If $\mathnormal{T}$ is injective then $\mathnormal{Y}$ is empty and the result trivially follows. Thus we may assume that $\mathnormal{T}$ is surjective.
To 
proceed with the proof we shall divide $\mathnormal{Y}$ into two 
pieces based on the following easy statement, and then deal with each piece separately.  Denote by $\D$ the contravariant functor $\Hom_k(\ ,k)$.

\begin{lemma}\label{L:isom} Let $V$ and $W$ be $k$-vector spaces and let $\alpha:V\otimes W\to\Hom_k(\D V,W)$ be the $k$-linear map given by $\alpha(v\otimes w)(f)=f(v)w, v\in V, w\in W, f\in\D V$ ($\alpha$ is an isomorphism if $\dim V<\infty$). For $x\in V\otimes W$ denote by $\End(V)x$ the cyclic $\End(V)$-submodule of $V\otimes W$ generated by $x$. Then $\End(V)x=V\otimes\Img\alpha(x)$.
\end{lemma}
\begin{proof} We have $x=\sum_{i=1}^s v_i\otimes w_i.$  If $s$ is the smallest possible, the sets of vectors $\{v_1,\dots,v_s\}$ and $\{w_1,\dots,w_s\}$ are linearly independent ~\cite[Theorem (1.2a), p. 142]{G1} so $\Img\alpha(x)$ is the span of $\{w_1,\dots,w_s\}$ and the rest is clear.
\end{proof}

\begin{definition} Let $\alpha_i:V_i\otimes W_i\to \Hom_k(\D V_i,W_i)$ be the $k$-linear map described in Lemma \ref{L:isom}, $i=1,\dots,l$.  If $x_i\in V_i\otimes W_i$ and $0\ne c\in k$, then $\Img\alpha_i(x_i)= \Img\alpha_i(cx_i)$, so if $p\in P$ is represented by $[x_1,\dots,x_l],\ x_i\in V_i\otimes W_i$, then for each $i $ the subspace $\Img\alpha_i(x_i)$  of $W_i$ is independent of the choice of representative for $p$. We set 
$\mathnormal{Y_1}=\{(\mathnormal{b},\mathnormal{p})\in\mathnormal{Y}:\sum_{i=1}^{l}(\dim{\mathnormal{V_i}})(\rank{\mathnormal{\alpha_i(x_i)}})<\dim{\mathnormal{U}}\}$ and 
$\mathnormal{Y_2}=Y-Y_1$.
\end{definition}

\begin{lemma}\label{L:union}$\mathnormal{Y}=\mathnormal{Y_1}\cup\mathnormal{Y_2}$ where
$\mathnormal{Y_1}$ is a closed subset of $\mathnormal{Y}$, and 
$\mathnormal{Y_2}$ is an open subset of 
$\mathnormal{Y}.$
\end{lemma}
\begin{proof} That 
$\mathnormal{Y}=\mathnormal{Y_1}\cup\mathnormal{Y_2}$ is obvious. For 
the other two statements consider the projection onto the second 
factor $\pi_2:\mathnormal{B}\times\mathnormal{P}\to\mathnormal{P}$ and note that
$\mathnormal{Y_1}$ ($\mathnormal{Y_2}$) is the intersection of 
$\mathnormal{Y}$ with the inverse image of the closed (open) subset 
of $\mathnormal{P}$ consisting of points corresponding to tuples of 
tensors satisfying
$\sum_{i=1}^{l}(\dim{\mathnormal{V_i}})(\rank{\mathnormal{\alpha_i(x_i)}})<(\geq)\dim{\mathnormal{U}}$.
\end{proof}
\begin{lemma}\label{L:dep} Assume that $\mathnormal{T}$ 
satisfies the right omnipresent maximal rank property. Suppose that 
$(\mathnormal{b},\mathnormal{p})\in\mathnormal{Y_1}$ and 
$\mathnormal{b}=[\varphi_1,\varphi_2,...,\varphi_l]$. Then for some 
$i$, $\varphi_i$ is not an isomorphism.
\end{lemma}
\begin{proof} If $x=[x_1,\dots,x_l]$ represents $p$, then $T(bx)=T\left([\varphi_1\cdot x_1,\dots,\varphi_l\cdot x_l]\right)=0$ because $(b,p)\in Y_1$.  By Lemma \ref{L:isom}, $Bx=\bigoplus_{i=1}^{\mathnormal{l}}\End(V_i)x_i=\bigoplus_{i=1}^{\mathnormal{l}}V_i\otimes\Img\alpha_i(x_i)$, so the restriction of $T$ to $Bx$ is injective because  $(b,p)\in Y_1$ and $T$ satisfies  the right omnipresent maximal rank property. Since $T(bx)=0$ then $bx=[\varphi_1\cdot x_1,\dots,\varphi_l\cdot x_l]=0$ whence $\varphi_i\cdot x_i=0$ for all $i$.
Because $\mathnormal{p}$ is a point in a projective space, at least 
one $x_i$ is not equal to 0. For this $i$, $\varphi_i$ 
is not an isomorphism. \end{proof}

\begin{lemma}\label{L:dim} Assume 
that $\mathnormal{T}$ satisfies the right omnipresent maximal rank 
property and that 
$\sum_{i=1}^{l}\dim\mathnormal{C(X_i)}\leq\dim\mathnormal{U}$. 
Suppose $\pi_2^{-1}(\mathnormal{X})\cap\mathnormal{Y_2}$ is nonempty. 
Then $\pi_2^{-1}(\mathnormal{X})\cap\mathnormal{Y_2}$ has Krull 
dimension at most $\Sigma(\dim{\mathnormal{V_i}})^2-1$, one less than 
the dimension of $\mathnormal{B}$.
\end{lemma}
\begin{proof} As with 
Lemma \ref{L:union} we consider the projection map onto the second 
factor $\pi_2:\mathnormal{B}\times\mathnormal{P}\to\mathnormal{P}$. 
Let $\mathnormal{X_2}\subset\mathnormal{X}$ be the set of points 
$\mathnormal{p}$ such that 
$\sum_{i=1}^{l}(\dim{\mathnormal{V_i}})(\rank{\alpha_i(x_i)})\geq\dim{\mathnormal{U}}$ for any $x=[x_1,\dots,x_l]$ representing $p$. Since
$\mathnormal{X_2}$ is open in $\mathnormal{X}$ and by assumption 
nonempty, $\dim{\mathnormal{X_2}}$ = $\dim{\mathnormal{X}}$. 
Pick any point $\mathnormal{p}$ in $\mathnormal{X_2}$ and, identifying $B$ with $B\times\{p\}$, consider 
the composite $\mathnormal{T':B\to U}$ of $T$ and the $k$-linear map $B\to \bigoplus_{i=1}^{\mathnormal{l}}(\mathnormal{V_i\otimes W_i})$ sending $b$ to $bx$. Clearly $\Ker \mathnormal{T'}=\pi_2^{-1}(\mathnormal{p})\cap\mathnormal{Y_2}$ and $\Img \mathnormal{T'}=T(Bx)=T(\bigoplus_{i=1}^{\mathnormal{l}}(\mathnormal{V_i\otimes \Img\alpha_i(x_i)}))$ (use Lemma \ref{L:isom}). 
 From the assumptions on the ranks of the $\alpha_i(x_i)$ and that 
$\mathnormal{T}$ satisfies the right omnipresent maximal rank property, 
we conclude that $\mathnormal{T'}$ is surjective. We then conclude 
that 
$\dim{(\pi_2^{-1}(\mathnormal{p})\cap\mathnormal{Y_2})}=\dim{\mathnormal{B}}-\dim{\mathnormal{U}}$.

Having 
computed the dimensions of the fibers of 
$\pi_2^{-1}(\mathnormal{X})\cap\mathnormal{Y_2}$ over 
$\mathnormal{X}_2$ we then see that the dimension of 
$\pi_2^{-1}(\mathnormal{X})\cap\mathnormal{Y_2}$ equals 
$\dim{\mathnormal{X}}+\dim{\mathnormal{B}}-\dim{\mathnormal{U}}$ = 
$\sum_{i=1}^{l}\dim\mathnormal{C(X_i)}-1+\dim{\mathnormal{B}}-\dim{\mathnormal{U}}\leq\dim{\mathnormal{B}}-1=\Sigma(\dim{\mathnormal{V_i}})^2-1$. 
\end{proof}
The proof of the theorem is now easy. 
\medskip\begin{flushleft}\emph{Proof of theorem.} By Lemma 
\ref{L:dep} the image of 
$\pi_2^{-1}(\mathnormal{X})\cap\mathnormal{Y_1}$ in $\mathnormal{B}$ 
will be contained in the proper closed subset of $\mathnormal{B}$ 
consisting of points where at least one $\varphi_i$ is not an isomorphism. By Lemma 
\ref{L:dim} $\pi_2^{-1}(\mathnormal{X})\cap\mathnormal{Y_2}$ has 
dimension less than that of $\mathnormal{B}$ and so its closure also 
does. Since $\pi_1$ is a projective morphism, the image in 
$\mathnormal{B}$ of the closure of 
$\pi_2^{-1}(\mathnormal{X})\cap\mathnormal{Y_2}$ will be closed and 
have dimension less than that of $\mathnormal{B}$.  Since by Lemma 
\ref{L:union} $\mathnormal{Y=Y_1\cup Y_2}$ we are done. 
\hskip1.04in$\square$
\end{flushleft} \medskip

\begin{corollary}\label{gen1} 
Assume that $\mathnormal{T}$ satisfies the right omnipresent maximal 
rank property. Make a choice of subspaces 
$\mathnormal{Z}_i\subset\mathnormal{V}_i\otimes\mathnormal{W}_i$, 
$i=1,...,l$. Then there exists a dense Zariski open subset $A\subset 
\mathnormal{B}$ such that if $[\varphi_1,...,\varphi_l]\in A$ then 
the 
restriction of $\mathnormal{T}$ to the direct sum of the 
$\varphi_i(\mathnormal{Z}_i)$ has maximal 
rank.
\end{corollary}
\begin{proof} We first do the case where 
$\sum_{i=1}^l\dim{(\mathnormal{Z}_i)}\le\dim{\mathnormal{(U)}}$. In 
Theorem \ref{T:main} set $\mathnormal{Z}_i=\mathnormal{C(X_i)}$. Choose $A$ to 
be the complement of any proper Zariski closed subset of 
$\mathnormal{B}$ containing 
$\pi_1(\pi_2^{-1}(\mathnormal{X})\cap\mathnormal{Y})$. For 
$[\varphi_1,...,\varphi_l]\in A$, 
$\bigoplus_{i=1}^l\varphi_i(\mathnormal{Z}_i)$ intersects the kernel 
of $\mathnormal{T}$ only in 0. Thus the restriction of 
$\mathnormal{T}$ to $\bigoplus_{i=1}^l\varphi_i(\mathnormal{Z}_i)$ is 
injective. When 
$\sum_{i=1}^l\dim{(\mathnormal{Z}_i)}=\dim{\mathnormal{(U)}}$ it is 
also surjective.

For the case where 
$\sum_{i=1}^l\dim{(\mathnormal{Z}_i)}>\dim{\mathnormal{(U)}}$ choose 
subspaces $Z'_i\subset\mathnormal{Z}_i$ such that 
$\sum_{i=1}^l\dim{(Z'_i)}=\dim{\mathnormal{(U)}}$. By 
the previous case we find $A$ such that if 
$[\varphi_1,...,\varphi_l]\in A$ then the restriction of 
$\mathnormal{T}$ to $\bigoplus_{i=1}^l\varphi_i(Z'_i)$ 
is surjective, so the restriction of $\mathnormal{T}$ to 
$\bigoplus_{i=1}^l\varphi_i(\mathnormal{Z}_i)$ is also surjective. 
\end{proof}

\begin{corollary}\label{gen2} Assume that 
$\mathnormal{T}$ is surjective and satisfies the right omnipresent maximal rank property. 
Fix integers $a_i,\ 0\le 
a_i\le(\dim{\mathnormal{V}_i})(\dim{\mathnormal{W}_i})$, $i=1,...,l$, 
such that $\sum a_i=\dim{U}$. For each $i$ choose $a_i$ linearly 
independent elements $m(i,j)$, $1\le j\le a_i$, of 
$\mathnormal{V}_i\otimes\mathnormal{W}_i$. Then there exists a dense 
Zariski open subset $A\subset \mathnormal{B}$ such that if 
$[\varphi_1,...,\varphi_l]\in A$ then the elements 
$\mathnormal{T}(\varphi_i(m(i,j)))$ form a basis for 
$\mathnormal{U}$.
\end{corollary}
\begin{proof} In Corollary 
\ref{gen1} set $\mathnormal{Z}_i$ equal to the span of the 
$m(i,j)$'s. 
\end{proof}

\begin{definition} \label{D:uni}
We say 
that $\mathnormal{T}$ satisfies the left general maximal rank 
property if and only if for a general choice of subspaces 
$\mathnormal{V'_i}\subset\mathnormal{V_i}$ for $i=1,...,l$ the 
restriction of $\mathnormal{T}$ to $\oplus_{i=1}^{l}(V'_i\otimes W_i)$ has 
maximal rank, that is, is either injective or 
surjective.
\end{definition}

By a general choice of subspaces we 
mean the following. Once the dimensions of the $\mathnormal{V'_i}$'s 
to be chosen are fixed, the set of all possible choices of 
$\mathnormal{V'_i}$'s can be identified with a product of 
Grassmanians. We mean that there exists a Zariski open dense subset 
of that product such that if the choice of  $\mathnormal{V'_i}$'s 
comes from that set, then the restriction of $\mathnormal{T}$ has 
maximal rank.

Similar to Definitions \ref{D1:uni} and \ref{D:uni} one can define what it means for the map $T$ of (\ref{mapT}) to satisfy the  left omnipresent or right general maximal rank property.  With these definitions, we leave it to the reader to interchange appropriately the words \lq\lq left" and \lq\lq right" in the above assertions and obtain true statements.  Of course, this comment also applies to the remainder of the section.

\begin{corollary}\label{rgen} If ${T}$ 
satisfies the right omnipresent maximal rank property then 
$\mathnormal{T}$ satisfies the left general maximal rank 
property.
\end{corollary}
\begin{proof} Make a choice of subspaces 
$\mathnormal{V'_i}\subset\mathnormal{V}_i$ for $i=1,...l$. In 
Corollary \ref{gen1} set 
$\mathnormal{Z}_i=\mathnormal{V'_i}\otimes\mathnormal{W}_i$. Notice 
that 
$\varphi_i(\mathnormal{V'_i}\otimes\mathnormal{W}_i)=\varphi_i(\mathnormal{V'_i})\otimes\mathnormal{W}_i$. 
A general tuple of endomorphisms $[\varphi_1,...,\varphi_l]$ applied 
to a specific tuple of subspaces 
$[\mathnormal{V'_1},...,\mathnormal{V'_l}]$ gives a general tuple of 
subspaces. 
\end{proof}

In Section 4 we will give examples to 
show that the right omnipresent maximal rank property does not imply 
the left omnipresent maximal rank property and the right general 
maximal rank property does not imply the left general maximal rank 
property.

We now indicate how to dualize the above results of this section. Let $V_1,V_2,...,V_l,W_1,W_2,...,W_l$, $Q$ be finite 
dimensional vector spaces and let 
\begin{equation}\label{mapS}\mathnormal{S}:Q\to\bigoplus_{i=1}^{\mathnormal{l}}(\mathnormal{V_i\otimes 
W_i})
\end{equation}
be a  linear transformation.
\begin{definition} 
\label{D:unidual}
We say that $\mathnormal{S}$ satisfies the right 
omnipresent maximal rank property if and only if for every choice of 
subspaces $\mathnormal{W'_i}\subset\mathnormal{W_i}$ for $i=1,...,l$ 
the composition of $S$ with the linear transformation 
$$\oplus_{i=1}^l(1_{V_i}\otimes\tau_i):\oplus_{i=1}^{\mathnormal{l}}(\mathnormal{V_i\otimes 
W_i})\to\oplus_{i=1}^{\mathnormal{l}}(\mathnormal{V_i\otimes 
(W_i/W'_i)})$$ 
has maximal rank, that is, is either injective or 
surjective, where $\tau_i:W_i\to W_i/W'_i$ is the natural projection.  And we say 
that $\mathnormal{S}$ satisfies the left general maximal rank 
property if and only if for a general choice of subspaces 
$\mathnormal{V'_i}\subset\mathnormal{V_i}$ for $i=1,...,l$ the composition of $S$ with the linear transformation 
$$\oplus_{i=1}^l(\sigma_i\otimes1_{W_i}):\oplus_{i=1}^{\mathnormal{l}}(\mathnormal{V_i\otimes 
W_i})\to\oplus_{i=1}^{\mathnormal{l}}(\mathnormal{(V_i/V'_i)\otimes 
W_i})$$ 
has maximal rank, where $\sigma_i:V_i\to V_i/V'_i$ is the natural projection.
\end{definition}

The following lemma shows that the question of whether a map of the type (\ref{mapT}) satisfies the omnipresent or general maximal rank property  is equivalent to the same question for a map of the type (\ref{mapS}).
\begin{lemma}\label{cross}  Let 

$$\begin{CD}
@.@.0\\
@.@.@VVV\\
@.@.B'\\
@.@.@ViVV\\
0@>>>A@>f>>B@>g>>C@>>>0\\
@.@.@VqVV\\
@.@.B''\\
@.@.@VVV\\
@.@.0
\end{CD}\ 
$$
be an exact diagram in an abelian category.  Then $gi$  is monic (epi) if and only if $qf$  is monic (epi).
\end{lemma}
\begin{proof} By the $3\times 3$ lemma the following commutative diagram is exact.

$$\begin{CD}
@.0@.0@.0\\
@.@VVV@VVV@VVV\\
0@>>>\Ker gi@>>>B'@>gi>>
\Img gi@>>>0\\
@.@VhVV@ViVV@VjVV\\
0@>>>A@>f>>B@>g>>C@>>>0\\
@.@VpVV@VqVV@VVV\\
0@>>>\Coker h@>r>>B''@>>>\Coker j@>>>0\\
@.@VVV@VVV@VVV\\
@.0@.0@.0
\end{CD} 
$$\vskip.1in
\noindent Hence $gi$ is monic if and only if $\Ker gi=0$, if and only if $p$ is iso, if and only if $qf$ is monic. The  rest of the proof is similar.
\end{proof}

\begin{corollary}
\label{kercoker}  $(\mathrm a)$ If a linear transformation $\mathnormal{T}:\bigoplus_{i=1}^{\mathnormal{l}}(\mathnormal{V_i\otimes 
W_i})\to\mathnormal{U}$ is surjective, it satisfies the left  omnipresent  (general) maximal rank property if and only if so does the inclusion $\Ker\mathnormal{T}\to \bigoplus_{i=1}^{\mathnormal{l}}(\mathnormal{V_i\otimes 
W_i})$.

$(\mathrm b)$ If a linear transformation $\mathnormal{S}:Q\to\bigoplus_{i=1}^{\mathnormal{l}}(\mathnormal{V_i\otimes 
W_i})$ is injective, then it satisfies the left  omnipresent (general) maximal rank property if and only if so does the projection $\bigoplus_{i=1}^{\mathnormal{l}}(\mathnormal{V_i\otimes 
W_i})\to\Coker S$.
\end{corollary}
\begin{proof}  The statement follows immediately from Lemma \ref{cross}. \end{proof}

We note that if the map $T$ of Corollary \ref{kercoker}(a) is injective, it satisfies both the right omnipresent and left general maximal rank property, and if $T$ is neither surjective nor injective then it satisfies neither of the properties.  A similar remark applies to the map $S$ of Corollary \ref{kercoker}(b).

A different way to relate the maps of the types (\ref{mapT}) and (\ref{mapS}) is through the vector space duality $\D$.

\begin{proposition}\label{dual} A linear transformation $T:\oplus_{i=1}^{l}(V_i\otimes W_i)\to U$ satisfies the left omnipresent (general) maximal  rank property if and only if so does its dual $\D T:\D U\to \oplus_{i=1}^{l}(\D V_i\otimes \D W_i)$.
\end{proposition}
\begin{proof} For $i=1,\dots,l$ let $X_i$ be a subspace of $V_i$ and $f_i:X_i\to V_i$ the inclusion map, then $\D f_i:\D V_i\to\D X_i$ is an epimorphism with $\Ker\D f_i=X_i^{\perp}=\{\phi\in\D V_i:\phi(X_i)=0\}$ so that $\D X_i\cong\D V_i/ X^{\perp}_i.$ Therefore $T\circ(\oplus_{i=1}^l(f_i\otimes1_{W_i}))$ is monic (epi) if and only if $(\oplus_{i=1}^l(\D f_i\otimes1_{\D W_i}))\circ\D T$ is epi (monic), if and only if $(\oplus_{i=1}^l(\psi_i\otimes1_{\D W_i}))\circ\D T$ is epi (monic), where $\psi_i:\D V_i\to\D V_i/X^{\perp}_i$  is the natural projection. Note that $X_i$ runs through the set of all subspaces of $V_i$ if and only if $X^{\perp}_i$ runs through the set of all subspaces of $\D V_i$.  Hence $T$ satisfies the left omnipresent maximal rank property if and only if so does $\D T$.  For a fixed sequence of nonnegative integers $d_i\le n_i=\dim V_i$, the $l$-tuple $(X_1,\dots,X_i,\dots,X_l)$ runs through a dense open set of the product of Grassmanians $\prod_{i=1}^lG(d_i,V_i)$ if and only if $(X^{\perp}_1,\dots,X^{\perp}_i,\dots,X^{\perp}_l)$ runs through the corresponding dense open set of the product of Grassmanians $\prod^l_{i=1}G(n_i-d_i,\D V_i)$ under the isomorphism that is the product of the natural isomorphisms $\D:G(d_i,V_i)\to G(n_i-d_i,\D V_i)$, see ~\cite[p. 200]{gh}.  Therefore $T$ satisfies the left general maximal rank property if and only if so does $\D T$. 
\end{proof}

We end this section with a 
lemma showing that the
right omnipresent maximal rank property puts a restriction on the relative 
sizes of the vector spaces involved.
\begin{lemma}\label{dimensions} $\mathrm{(a)}$ If 
$\mathnormal{T}:\oplus^l_{i=1}V_i\otimes W_i\to U$ satisfies 
the right omnipresent maximal rank property and ${T}$ is surjective but not 
injective, then 
$\dim{{U}}<\sum^l_{i=1}(\dim{{V_i}})^2$.

$\mathrm{(b)}$  If 
${S}:Q\to\oplus^l_{i=1}V_i\otimes W_i$ satisfies 
the right omnipresent maximal rank property and ${S}$ is injective but not 
surjective, then 
$\dim{{Q}}<\sum^l_{i=1}(\dim{{V_i}})^2$.
\end{lemma}
\begin{proof} (a)
Suppose to the contrary that 
$\dim{{U}}\ge\sum^l_{i=1}(\dim{{V_i}})^2$. Let 
$\{v_{i1},...,v_{ik_i}\}$ be a basis for ${V_i}$. Express some 
nonzero element of $\Ker\mathnormal{T}$ in the form 
\newline\centerline{$\left(\sum_{j=1}^{k_1} v_{1j}\otimes w_{1j},\dots,\sum_{j=1}^{k_i} v_{ij}\otimes w_{ij}\dots,\sum_{j=1}^{k_l} v_{lj}\otimes w_{lj}\right)$} 
and let ${W'_i}$ equal the span 
of $\{w_{i1},...,w_{ik_i}\}$.  Then 
$\dim{\oplus^l_{i=1}({V_i}\otimes\mathnormal{W'_i})}\le\dim{U}$ but the 
restriction of ${T}$ to 
$\oplus^l_{i=1}({V_i}\otimes{W'_i})$ is neither surjective nor injective because its kernel is not zero.

(b) Follows from (a) and Proposition \ref{dual}.
\end{proof}

\section{Almost Split 
Morphisms and Preprojective Algebras}

We apply the results of Section 2 to representations of algebras which  provide a large supply of linear transformations of the form 
$\oplus_{i=1}^{l}(V_i\otimes W_i)\to U$ or $Q\to\oplus_{i=1}^{l}(V_i\otimes W_i)$.  Let $\Lambda$ be an associative $k$-algebra, let $\module\La$ be the category of finite dimensional left $\La$-modules, and let $g:B\to C$ and $f:A\to B$ be morphisms in  $\module\La$.  Replacing $B$ with an isomorphic module if necessary, we may assume that $B=V_1^{n_1}\oplus\dots\oplus V_l^{n_l}$ where $V_1,\dots,V_l$ are nonisomorphic indecomposable $\La$-modules, $l, n_1,\dots,n_l$ are nonegative integers, and $V^m$ stands for the direct sum of $m$ copies of $V$.  For $i=1,\dots,l$ denote by $W_i$ the $k$-space with a basis $e_{i1},\dots,e_{in_i}$, and for each $j=1,\dots,n_i$ denote by $h_{ij}:V_i\to V_i\otimes ke_{ij}$ the isomorphism of $\La$-modules sending each $v\in V_i$ to $v\otimes e_{ij}$.  Let \begin{equation}\label{h} h:B\to \oplus_{i=1}^{l}(V_i\otimes W_i)\end{equation} 
be the  isomorphism in $\module\La$ induced by the $h_{ij}$'s.  Denote by $g_{ij}:V_i\to C$ and $f_{ij}:A\to V_i$ the morphisms in $\module\La$ induced by $g$ and $f$, respectively, and consider the morphisms $T_i:V_i\otimes W_i\to C$ and $S_i:A\to V_i\otimes W_i$ defined by $T_i(v\otimes e_{ij})=g_{ij}(v), v\in V_i,$ and $S_i(a)=(f_{ij}(a)\otimes e_{ij}), a\in A,$ respectively.  Let 
\begin{equation}\label{maps} T:\oplus_{i=1}^{l}(V_i\otimes W_i)\to C\mathrm{\phantom{X} and \phantom{X}}S:A\to\oplus_{i=1}^{l}(V_i\otimes W_i)\end{equation} 
be the  morphisms in $\module\La$ induced by the $T_i$'s and $S_i$'s, respectively.  It is straight forward to check that 
\begin{equation}\label{formulae}g=Th\mathrm{\phantom{X} and \phantom{X}}S=hf.\end{equation}
\begin{proposition}\label{L:Kleiner}
Let $g:B\to C$ and $f:A\to B$ be morphisms in $\module\La$ with $B=\oplus^l_{i=1}V_i^{n_i}$ where $V_1,\dots,V_l$ are nonisomorphic indecomposable $\La$-modules. Let $h$ be the isomorphism in (\ref{h}), let $T$ and $S$ be the morphisms in (\ref{maps}) constructed from $g$ and $f$, respectively.\vskip.06in

If $g$ is a minimal right almost split morphism in $\module\La$ then:
\begin{itemize}
\item[$\mathrm{(a)}$] $T$ satisfies the right omnipresent maximal rank property.
\item[$\mathrm{(b)}$] For a general choice of $k$-subspaces $U_i\subset V_i$, the restriction of $g$ to $\oplus^l_{i=1}U_i^{n_i}$ has maximal rank.
\item[$\mathrm{(c)}$] If $g$ is surjective then $\dim C<\sum^l_{i=1}(\dim V_i)^2$.
\end{itemize}

If $f$ is a minimal left almost split morphism in $\module\La$ then:
\begin{itemize}
\item[$\mathrm{(d)}$]  $S$ satisfies the right omnipresent maximal rank property.
\item[$\mathrm{(e)}$] For a general choice of $k$-subspaces $U_i\subset V_i$, denote by $\sigma_i:V_i\to V_i/U_i$ the natural projection.  Then the  linear transformation $(\oplus^l_{i=1}\sigma_i^{n_i})\circ f$ has maximal rank.
\item[$\mathrm{(f)}$] If $f$ is injective then $\dim A<\sum^l_{i=1}(\dim V_i)^2$.
\end{itemize}

If $\,0\to A\overset{f}\to B\overset{g}\to C\to0$ is an almost split sequence in $\module\La$ then:
\begin{itemize}
\item[$\mathrm{(g)}$] $\dim B<2\sum^l_{i=1}(\dim V_i)^2-1$.
\end{itemize}
\end{proposition}
\begin{proof} (a) \hskip.03in Since $g$ is minimal right almost split, so is $T$ by (\ref{formulae}).  If $W_i'$ is a subspace of $W_i$, the $\La$-module $V_i\otimes W_i'$ is a direct summand of $\mathnormal{V_i\otimes W_i}$. Hence the restriction of $T$ to $\oplus_{i=1}^{\mathnormal{l}}(V_i\otimes W_i')$ is an irreducible morphism and thus  is either a monomorphism or an epimorphism ~\cite[Ch. V, Theorem 5.3(a) and Lemma 5.1(a)]{ars}, so (a) holds.  According to Corollary \ref{rgen}, $T$ satisfies the left general maximal rank property.   In view of the structure of the isomorphisms $h_{ij}$ constructed above, we conclude that (b) holds. Part  (c) is a direct  consequence of (a), formula (\ref{formulae}), and Lemma \ref{dimensions}(a).

(d)  The proof is similar to that of (a) using the analogous properties of minimal left almost split morphisms.

(e) If $f$ is surjective, the statement is clear.  If $f$ is not surjective, it is injective, and so is $S$ in view of formulas (\ref{formulae}).   By (d) and Corollary \ref{kercoker}(b), the projection $\oplus_{i=1}^{l}(V_i\otimes W_i)\to\Coker S$ satisfies the right omnipresent maximal rank property. By Corollary \ref{rgen}, it satisfies the left general maximal rank property, and so does $S$ by Corollary \ref{kercoker}(b). Then $f$ satisfies the desired property in view of formulas (\ref{formulae}). 

Another way to prove (d) and (e) is to note that both $\D f$ and $\D S$ are minimal right almost split morphisms in $\module\La^{\mathrm{op}}$, and then use (a), Corollary \ref{rgen}, and Proposition \ref{dual} together with formulas (\ref{formulae}).

(f) \hskip.03in The proof is similar to that of (c), using Lemma \ref{dimensions}(b).

(g) \hskip.03in The formula follows from (c) and (f).
\end{proof}

\begin{rmk}\label{irreducible} (a) Lemma \ref{dimensions} holds when  $k$ is an arbitrary field.  Hence so do parts (a), (c), (d), (f), and (g) of Proposition \ref{L:Kleiner}; moreover, they hold if 
$\module\La$ is replaced by any full subcategory of an abelian category closed under  extensions and direct summands where the objects and morphism sets are finite dimensional $k$-vector spaces and composition of morphisms is $k$-bilinear.

(b) Parts (a), (b), and (c)  of Proposition \ref{L:Kleiner} hold if  $g:B\to C$  is an irreducible morphism with $C$ indecomposable, and parts (d), (e), and (f) hold if $f:A\to B$ is an irreducible morphism with $A$ indecomposable.  This follows from the observation after Definition \ref{D1:uni} that  the right omnipresent maximal rank property of a  linear transformation is inherited by its appropriate restrictions, and from the dual statement.

(c) Parts (c), (f), and (g) of Proposition \ref{L:Kleiner} imply that for a fixed number of nonisomorphic indecomposable summands of the middle term of an almost split sequence, the summands cannot be much smaller than the end terms of the sequence, i.e., the multiplicities of the summands cannot be too large, and that there is a balance between the sizes of the end terms.  Part (c) is false if the morphism $g$ is not surjective, and part (f) is false if the morphism $f$ is not injective.
\end{rmk}

We will apply this in particular to the preprojective algebra where the grading introduced in \cite{K1} allows us to interpret the multiplication-by-arrow maps into (from) a fixed homogeneous component as a minimal right (left) almost split morphism of modules over the path algebra of the quiver.   We recall some facts from the latter paper.

For the remainder of this paper we fix a finite quiver 
$\Gamma=(\Gamma_{0},\Gamma_{1})$ 
 without oriented cycles with 
the set of vertices $\Gamma_{0}$ and the set of 
 arrows 
$\Gamma_{1}$. Let $\bar\Gamma=(\bar\Gamma_{0},\bar\Gamma_{1})$ be a new quiver 
 with 
$\bar\Gamma_{0}=\Gamma_{0}$
 and 
$\bar\Gamma_{1}=\Gamma_{1}\cup\Gamma_{1}^{*}$,
where 
$\Gamma_{1}\cap\Gamma_{1}^{*}=\emptyset$ and the elements of 
$\Gamma_{1}^{*}$ are in the following one-to-one 
correspondence 
with the elements of $\Gamma_{1}$: for 
each $\gamma:t\to v$ in $ 
\Gamma_{1}$, there is a unique element 
$\gamma^{*}:v\to t$ in 
$\Gamma_{1}^{*}$.  To turn the path 
algebra  $k\bar\Gamma$ of 
$\bar\Gamma$ over a field $k$  
into a graded $k$-algebra, we assign 
degree  0 to each trivial 
path $e_{t}, t\in\Gamma_{0}$, and each 
arrow $\gamma\in\Gamma_{1}$; 
degree 1 to each arrow 
$\gamma^{*}\in\Gamma_{1}^{*}$; and compute 
the degree of a 
nontrivial path $q=\delta_{1}\dots\delta_{r}$ as 
deg 
$q=\sum_{i=1}^{r}\text{deg}\,\delta_{i}$. Clearly, $k\Gamma$ is 
the $k$-subalgebra of $k\bar\Gamma$ comprising the elements of degree 0.  

Let $\mathbb N$ be the set of nonnegative integers. 
For all 
$t\in\Gamma_{0},\ d\in\Bbb N$, let $W^{t}_{d}$ be
the span of all those paths in $\bar\Gamma$ of degree $d$ that start 
at 
$t$. Note that  $W^{t}_{d}\in\module k\Gamma$ so
\begin{equation}\label{dirsum}
k\bar\Gamma=\underset{d\in\Bbb 
N}\oplus\underset{t\in\Gamma_{0}}\oplus 
W^{t}_{d}
\end{equation}
is a decomposition of $k\bar\Gamma$ as a direct sum of its 
left $k\Gamma$-submodules

Let now $a$ and $b$ be any two
functions $\Gamma_{1}\to k$ satisfying $a(\gamma)\ne0$ and 
$b(\gamma)\ne0$ for all $\gamma\in\Gamma_{1}$.  If  $s(\gamma)$ is the starting point and $e(\gamma)$ is the end point of $\gamma\in\Gamma_{1}$, for each $t\in\Gamma_{0}$ set 
$$
m_{t}=\underset{s(\gamma)=t}
{\underset{\gamma\in\Gamma_{1}}\sum}a(\gamma)\gamma^{*}\gamma-
\underset{e(\gamma)=t}{\underset{\gamma\in\Gamma_{1}}\sum} 
b(\gamma)\gamma\gamma^{*}
$$ 
and denote by $J$ the two-sided ideal 
of 
$k\bar\Gamma$ generated by the 
element
$$
\underset{t\in\Gamma_{0}}\sum m_{t}=\underset{\gamma\in\Gamma_{1}}\sum 
[\gamma^{*},\gamma]_{a,b}
$$ 
where $[\gamma^{*},\gamma]_{a,b}=a(\gamma)\gamma^{*}\gamma-b(\gamma)\gamma
\gamma^{*}*$ is the $(a,b)$-commutator of $\gamma^{*}$ and $\gamma$. 
The factor 
algebra 
${\mathcal P}_{k}(\Gamma)_{a,b}=k\bar\Gamma/J$ 
is the 
$(a,b)$-preprojective algebra of $\Gamma$. 

Since the elements $m_{t}$ are homogeneous 
of degree 1, $J$ is a homogeneous 
ideal containing no nonzero 
elements of degree 0.  Hence ${\mathcal 
P}_{k}(\Gamma)_{a,b}$ 
is 
a graded $k$-algebra, and the restriction to $k\Gamma$ of the natural 
projection 
$\pi:k\bar\Gamma\to{\mathcal P}_{k}(\Gamma)_{a,b}$ is 
an isomorphism of $k\Gamma$ 
with the
 subalgebra of ${\mathcal 
P}_{k}(\Gamma)_{a,b}$ comprising the elements of degree 0;
  we 
view the isomorphism as identification. From (\ref{dirsum}) we get 
$$
{\mathcal 
P}_{k}(\Gamma)_{a,b}=\underset{d\in\mathbb 
N}\oplus\underset{t\in\Gamma_{0}}\oplus 
V^{t}_{d}
$$
where $V^{t}_{d}=\pi(W^{t}_{d})\in\module k\Gamma$.
If 
$\gamma\in\Gamma_{1}$ we write $\beta=\pi(\gamma)$ and 
$\beta^{*}=\pi(\gamma^{*})$.  If $q$ is a path 
in $\bar\Gamma$ 
starting at $t$ and ending at $v$, we 
call $\pi(q)$ a path in 
${\mathcal P}_{k}(\Gamma)_{a,b}$ starting at $t$ and ending
 at 
$v$.
 Then  $V^{t}_{d}$ is the span of all paths of degree $d$ in 
 ${\mathcal P}_{k}(\Gamma)_{a,b}$ starting at $t$.  Since we 
identify $k\Gamma$ 
 with $\pi(k\Gamma)$, we in particular identify 
$e_{t}$ with 
 $\pi(e_{t}),\ t\in\Gamma_{0}$; $\gamma$ with 
$\beta=\pi(\gamma),\ 
 \gamma\in\Gamma_{1}$; $W^{t}_{0}$ 
with $V^{t}_{0}$; and  we set $W^{\,t}_{-1}=V^{\,t}_{-1}=0$.

We need the following statement. When appropriate, 
the map $(c):X\to Y$ denotes the right multiplication by $c$.
\begin{theorem}\label{negref3} Suppose $V^{t}_{d}\ne 0$ where $t\in\Gamma_0$, $d\in{\mathbb N}$. 
\begin{itemize}

\item[$\mathrm{(a)}$] $V^{t}_{d}$ is indecomposable in $\module k\Gamma$, and  $V^{t}_{d}\cong V^{s}_{c}$ in $\module k\Gamma$, $s\in\Gamma_0$, $c\in{\mathbb N}$,  if and only if $t=s$ and $d=c$.

\item[$\mathrm{(b)}$] The map
\hskip.05in${g^{t}_{d}}:\left(\underset{s(\gamma)=t}{\oplus} 
V^{e(\gamma)}_{d}\right)\oplus\left(\underset{e(\gamma)=t}\oplus 
V^{s(\gamma)}_{d-1}\right)\longrightarrow V^{t}_{d}$\hskip.07in 
induced by the right multiplications
$(a(\gamma)\beta):V^{e(\gamma)}_{d}\to V^{t}_{d}$, $s(\gamma)=t$, and 
$(-b(\gamma)\beta^{*}):V^{s(\gamma)}_{d-1}\to V^{t}_{d}$,  $e(\gamma)=t$, where ${\gamma
\in\Gamma_{1}}$, is a minimal right almost split morphism in $\module k\Gamma$.

\item[$\mathrm{(c)}$] The map
\hskip.05in$f^{t}_{d}: V^{t}_{d}\longrightarrow 
\left(\underset{s(\gamma)=t}{\oplus} 
V^{e(\gamma)}_{d+1}\right)\oplus\left(\underset{e(\gamma)=t}{\oplus }
V^{s(\gamma)}_{d}\right)$\hskip.05in induced by the right multiplications 
$(\beta^{*}¥):V^{t}_{d}\to V^{e(\gamma)}_{d+1}¥$, $s(\gamma)=t$, and  
$(\beta):V^{t}_{d}\to V^{s(\gamma)}_{d}¥$,  
$e(\gamma)=t$, where $\gamma\in\Gamma_{1}¥$, is a  minimal left almost split morphism in $\module k\Gamma$.

\item[$\mathrm{(d)}$] If $V^{t}_{d+1}\ne 0$ then\hskip.1in $0\to V^{t}_{d}\overset{f^{t}_{d}}\longrightarrow 
\left(\underset{s(\gamma)=t}{\oplus} 
V^{e(\gamma)}_{d+1}\right)\oplus\left(\underset{e(\gamma)=t}{\oplus }
V^{s(\gamma)}_{d}\right)\overset{g^{t}_{d+1}}\longrightarrow  V^{t}_{d+1}\to0$ is an almost split sequence in $\module k\Gamma$.
\end{itemize}
\end{theorem}
\begin{proof} These are parts of ~\cite[Theorem 1.1 and Corollary 1.3]{K1} combined with well-known properties of preprojective modules, see ~\cite[VIII.1]{ars}.
\end{proof}

Applying parts (b) and (d) of Proposition \ref{L:Kleiner} to Theorem \ref{negref3}, we obtain the following statement.
\begin{corollary}\label{leftgeneral}
$\mathrm{(a)}$ In the setting of Theorem \ref{negref3}(b), for a general choice of $k$-subspaces
\hskip.05in $U^{e(\gamma)}_d\subset V^{e(\gamma)}_d$ and \hskip.05in $U^{s(\gamma)}_{d-1}\subset V^{s(\gamma)}_{d-1}$, the restriction of \hskip.02in$g^t_d$\hskip.05in to $(\underset{s(\gamma)=t}{\oplus} 
U^{e(\gamma)}_{d})\oplus(\underset{e(\gamma)=t}\oplus 
U^{s(\gamma)}_{d-1})$ has maximal rank.

$\mathrm{(b)}$ In the setting of Theorem \ref{negref3}(c), for a general choice of $k$-subspaces
\hskip.05in $U^{e(\gamma)}_{d+1}\subset V^{e(\gamma)}_{d+1}$ and \hskip.05in $U^{s(\gamma)}_{d}\subset V^{s(\gamma)}_{d}$, denote by $\sigma^{e(\gamma)}_{d+1}:V^{e(\gamma)}_{d+1}\to V^{e(\gamma)}_{d+1}/U^{e(\gamma)}_{d+1}$ and $\sigma^{s(\gamma)}_{d}:V^{s(\gamma)}_{d}\to V^{s(\gamma)}_{d}/U^{s(\gamma)}_{d}$ the natural projections.  Then the  linear transformation \hskip.05in$((\underset{s(\gamma)=t}{\oplus} 
\sigma^{e(\gamma)}_{d+1})\oplus(\underset{e(\gamma)=t}{\oplus }
\sigma^{s(\gamma)}_{d}))\circ f^t_d$\hskip.07in has maximal rank.
\end{corollary}

\begin{rmk} As follows from Remark \ref{irreducible}(b), if one leaves out any number of summands in the direct sum of part (b) of Theorem \ref{negref3} and replaces the map $g^t_d$ by its restriction to the sum of the remaining summands, Corollary \ref{leftgeneral}(a) will still hold.  Likewise, if one leaves out any number of summands in the direct sum of part (c) of Theorem \ref{negref3} and replaces the map $f^t_d$ by its composition with the projection onto the sum of the remaining summands, Corollary \ref{leftgeneral}(b) will still hold.
\end{rmk}

The results of this section have dealt with left modules over a $k$-algebra $\La$ and with the right multiplication-by-arrow maps in the preprojective algebra. One may ask if analogous results are true for right $\La$-modules and for the left mulitplication-by-arrow maps.  We leave it to the reader to state the analog of Proposition \ref{L:Kleiner}, and note that  ~\cite[Theorem 1.1 and  Corollary 1.3]{K1} address  left multiplication by arrows in ${\mathcal P}_{k}(\Gamma)_{a,b}$ by replacing $W^t_d$ and $V^t_d$ with $W_{t,d}$, the span of all those paths in $\bar\Gamma$ of degree $d$ that end at $t$, and $V_{t,d}=\pi(W_{t,d})$, respectively.  Since $V_{t,d}$ is a finite dimensional right $k\Gamma$-module for all $t$ and $d$, with the appropriate replacements the analogs of Theorem \ref{negref3} and Corollary \ref{leftgeneral} hold.  These remarks also apply to the considerations of Section 4.

\section{Corollaries and Examples}
In this section we strengthen Corollary \ref{leftgeneral} in a form that is analogous to the result of Hochster and Laksov \cite{HL1}.  To help the reader see the analogy we shall first state their result.

Set ${R=k[x_1,x_2,...,x_r]}$, the commutative 
polynomial ring graded by degree, and denote by ${R_d}$ its 
 homogeneous piece  of degree $\mathnormal{d}$. Let 
$\mathnormal{N(r,d)}$ be the dimension of $\mathnormal{R_d}$ as a 
vector space over $\mathnormal{k}$. The following is then the result 
of Hochster and Laksov \cite{HL1}.
\begin{theorem}\label{HL}
Given an 
integer ${d\geq2}$, we determine an integer 
$\mathnormal{n}$ by the 
inequalities
\[{(n-1)r<N(r,d+1)\leq nr}\]
and let 
${s=N(r,d+1)-(n-1)r}$. Then if 
${F_1,F_2,...,F_n}$ are ${n}$ general forms in 
${R_d}$ we have that the ${(n-1)r}$ forms 
${x_jF_i}$ for ${j=1,...,r}$ and 
${i=1,2,...,n-1}$ together with the ${s}$ forms 
${x_jF_n}$ for ${j=1,2,...s}$ (in total 
${N(r,d+1)}$ forms) are a ${k}$-vector space 
basis for ${R_{d+1}}$.
\end{theorem}
By ``general forms" 
they mean that there exists a dense Zariski open subset of the affine 
space $\mathnormal{(R_d)^n}$ such that if the $\mathnormal{n}$-tuple 
$\mathnormal{(F_1,F_2,...,F_n)}$ is chosen from that open set, then 
the conclusion follows.

We wish to apply Corollary 
\ref{gen2} to the maps $g^t_d$ of Theorem \ref{negref3}, which is possible according to Proposition \ref{L:Kleiner}(a). To make the result clearly analogous to the result of Hochster and Laksov we must set up our notation properly.

Fix a vertex 
$t\in\Gamma_0$ and  a nonnegative integer $d$. Let $s_1,s_2,...,s_m$ be the distinct vertices that 
have in $\Gamma_1$ arrows from them to $t$, and let $u_1,u_2,...u_n$ be 
the distinct vertices with arrows in $\Gamma_1$ going from $t$ to 
them.  To match things up with the 
set up in Section 2, for $i=1,2,...,m$ let 
${V_i}=V_{d-1}^{s_i}$, let ${W_i}$ be the 
$k$-linear span of the arrows $\beta^*_{i,j}$ in $\Gamma_1^*$ going from $t$ to $s_i$, and set $w_{i,j}=-b(\beta_{i,j})\beta^*_{i,j}$.  For $i=m+1,m+2,...,m+n=l$ let 
${V_i}=V_{d+1}^{u_{i-m}}$, let ${W_i}$ be 
the $k$-linear span of the arrows $\beta_{i,j}$ in $\Gamma_1$ going from $t$ to 
$u_{i-m}$, and set $w_{i,j}=a(\beta_{i,j})\beta_{i,j}$. For all $i$, we choose $\{w_{i,j}\}$
as a  basis for ${W_i}$ and put the  $w_{i,j}$'s in a column vector $x_i$.  Set ${U}=V_{d}^t$. 
Let ${M_i}$ be the vector space of 
$\dim{{V_i}}\times\dim{{W_i}}$ matrices with 
elements in $k$. Let $B'$ be the affine space $\prod_{i=1}^l{V}_i^{\dim{{V}_i}}$. An element 
${b'}$ of $B'$ is an $l$-tuple $[{b'_1},{b'_2},...,{b'_l}]$ 
where each $b'_i$ is a $\dim{{V}_i}$-tuple of elements of ${V}_i$, written as a row vector. For $b'\in B'$ and $m(i)\in M_i$, using ordinary matrix multiplication and the multiplication and addition in the preprojective algebra,  we see that ${b'_im(i)x_i}$  is an element of ${U}=V_{d}^t$. 
\begin{corollary}\label{AHL} Let $d>0$ and $V^t_d\ne0$.
Fix integers $a_i$ satisfying $\,0\leq 
a_i\leq(\dim{\mathnormal{V_i}})(\dim{\mathnormal{W_i}})$, 
$i=1,...,l$, and $\sum a_i$ = $\dim{\mathnormal{U}}$. For 
each $i$ choose $a_i$ linearly independent elements of  ${M_i}$ and call them $m(i,j)$, $1\leq j \leq 
a_i$.
There exists a Zariski open 
dense subset ${E}$ of ${B'}$ such that if 
${b'}=[{b'_1},{b'_2},...,{b'_l}]\in{E}$, 
then the elements ${b'_im(i,j)x_i}$, $1\leq i\leq l$, $1\leq 
j\leq a_i$, form a basis for 
$\mathnormal{U}=V_{d}^t$.
\end{corollary}
\begin{proof} We already have a chosen basis for each ${W}_i$. Suppose we also choose a basis for each ${V}_i$. The pairwise tensor products of these basis elements give a basis for ${V}_i\otimes{W}_i$, so we may identify  ${V}_i\otimes{W}_i$ with $M_i$. We may also identify $\End{({V}_i)}$ with ${V}_i^{\dim{{V}_i}}$ by matching $\varphi_i\in \End{({V}_i)}$ with the image under $\varphi_i$ of the chosen basis. Under these identifications the elements ${T}(\varphi_i(m(i,j)))$ appearing in Corollary \ref{gen2} become identified with the elements $\mathnormal{b'_im(i,j)x_i}$ appearing in Corollary \ref{AHL}. Thus Corollary \ref{AHL} is a particular case of Corollary \ref{gen2}. 
\end{proof}

With the proper choice of the $m(i,j)$ we can get a 
corollary that sounds even more like the result of Hochster and 
Laksov.
\begin{corollary}\label{AAHL} Let $d>0$ and $V^t_d\ne0$.
For each $i$ satisfying $1\leq i\leq m$, let 
$\beta_{i,j}^*$, $1\leq j\leq\dim{\mathnormal{W}_i}$, be the new 
arrows going from $t$ to $s_i$. For each $i$ satisfying $m+1\leq i\leq l$, let 
$\beta_{i,j}$, $1\leq j\leq\dim{\mathnormal{W}_i}$, be the old arrows 
going from $t$ to $u_{i-m}$. Choose positive integers $n_i$, $1\leq 
i\leq l$, satisfying $1\le n_i\le\dim V_i$ and
\newline\centerline{$\sum_{i=1}^l(n_i-1)\dim{{W}_i}<\dim{V}_{d}^t\leq\sum_{i=1}^ln_i\dim{\mathnormal{W}_i},$}\medskip

\noindent and set 
$c$ = 
$\dim{V}_{d}^t-\sum_{i=1}^l(n_i-1)\dim{\mathnormal{W}_i}$. 
Write $c$ as a sum of nonnegative integers $c=c_1+c_2+...+c_l$, 
$0\leq c_i\leq\dim{{W}_i}$. For a general choice of
$\sum_{i=1}^ln_i$ elements  $F_{i,k}$, $1\le i\le l$, $1\le k\le n_i$, where $F_{i,k}\in {V}_{d-1}^{s_i}$ for $1\leq i\leq m$ and $F_{i,k}\in{V}_{d}^{u_{i-m}}$ for $m+1\leq i\leq l$, the following
$\dim{{V}_{d}^t}$ elements 
form a basis for ${V}_{d}^t$:\medskip
\newline $F_{i,k}\beta_{i,j}^*$ for $1\leq i\leq m$,
$1\leq k \leq n_i-1$, $1\le j\le\dim{\mathnormal{W}_i}$;
\newline $F_{i,n_i}\beta_{i,j}^*$ for $1\leq i\leq m$, $1\leq j \leq
c_i$;
\newline $F_{i,k}\beta_{i,j}$ for $m+1\leq i\leq l$, $1\leq k \leq n_i-1$, $1\le j\le\dim{\mathnormal{W}_i}$;  
\newline $F_{i,n_i}\beta_{i,j}$ for $m+1\leq i\leq l$, $1\leq j \leq c_i$.\medskip
\newline Here in an inequality giving
the range of possible $j$ or $k$, if the number on the right is less than 1,
we simply mean there are no such $j$ or $k$.  
\end{corollary} 
\begin{proof} 
Choose the $m(i,j)$'s as follows. Note that 
$a_i=(n_i-1)\dim{\mathnormal{W}_i}+c_i$. For a fixed $i$, the $a_i$ elements
$m(i,j)$ will be the $(n_i-1)\dim{{W}_i}$ distinct 
matrices having a 1 in one place among the 
$(n_i-1)\dim{{W}_i}$ positions available in the first 
$(n_i-1)$ rows of the 
$\dim{{V}_i}\times\dim{{W}_i}$ matrices 
involved and zeros elsewhere. The remaining $c_i$ elements $m(i,j)$ have a 1 
in one of the first $c_i$ places in the $n_i$-th rows, and zeros 
elsewhere.
\end{proof}
\begin{corollary} $\mathrm{(a)}$ If $d>0$ and $V^t_d\ne0$ then $\dim V^t_d<\sum^n_{j=1}\left(\dim V_d^{u_j}\right)^2+\sum^m_{i=1}\left(\dim V_{d-1}^{s_i}\right)^2$.\vskip.03in

\noindent If $d\ge0$ and $V^t_{d+1}\ne0$ then:\vskip.04in

$\mathrm{(b)}$  $0<\dim V^t_d<\sum^n_{j=1}\left(\dim V_{d+1}^{u_j}\right)^2+\sum^m_{i=1}\left(\dim V_{d}^{s_i}\right)^2$.\vskip.03in

$\mathrm{(c)}$ $0<\left(\underset{s(\gamma)=t}{\sum} 
\dim V^{e(\gamma)}_{d+1}\right)+\left(\underset{e(\gamma)=t}\sum 
\dim V^{s(\gamma)}_{d}\right)<2\left(\sum^n_{j=1}\left(\dim V_{d+1}^{u_j}\right)^2+\sum^m_{i=1}\left(\dim V_{d}^{s_i}\right)^2\right)-1$\vskip.05in \noindent where $\gamma\in\Gamma_1$. 
\end{corollary}
\begin{proof}  This is a direct consequence of Theorem \ref{negref3} and parts (c), (f), and (g) of Proposition \ref{L:Kleiner}.
\end{proof}

\begin{example} This example
shows that the $F_{i,j}$ of Corollary \ref{AAHL} must be chosen generically. In other words the right omnipresent maximal rank property does not imply the left omnipresent maximal rank property. Let the quiver 
$\Gamma$ have two vertices labeled 1 and 2 and one arrow $\beta$ 
going from 1 to 2. $\bar{\Gamma}$ then has in addition one new arrow 
$\beta^*$ going from 2 to 1. For any choice of nonzero functions $a$ 
and $b$ the relations become $\beta\beta^*=\beta^*\beta=0$. In Theorem \ref{negref3} set $d=1$ and $t=2$. The map becomes
$V_0^1\to V_1^2$ where
$V_0^1$ has basis $\{e_1,\beta\}$, and $V_1^2$ has basis 
$\{\beta^*\}$. The map is multiplication by $\beta^*$ so $e_1$ 
goes to $\beta^*$ and $\beta$ goes to 0. Consider one dimensional 
subspaces of $V_0^1$. The one spanned by $\beta$ maps to 0 and so 
does not surject onto $V_1^2$, all others do surject onto 
$V_1^2$.
\end{example}

\begin{example} Here we show that if in Theorem 
\ref{T:main} the hypothesis that $\mathnormal{T}$ satisfies the right 
omnipresent maximal rank property is weakened to the right general maximal rank property, then the conclusion might 
not follow. In other words the right general maximal rank property does not imply the left general maximal rank property. Let $\mathnormal{V}$ be a vector space of dimension 3 
with basis $\{v_1,v_2,v_3\}$. Let $\mathnormal{W}$ be a vector space 
of dimension 2 with basis $\{w_1,w_2\}$. Let $\mathnormal{U}$ be the 
quotient of $\mathnormal{V\otimes W}$ by the subspace spanned by 
$\{v_1\otimes w_1,v_2\otimes w_1\}$. Finally let 
$\mathnormal{T}:\mathnormal{V\otimes W}\to\mathnormal{U}$ be the 
quotient map. The only one-dimensional subspace $\mathnormal{W}'$ of $\mathnormal{W}$ such that $\mathnormal{V}\otimes\mathnormal{W}'$ has nonzero intersection with the kernel of $\mathnormal{T}$ is the span of $w_1$. Thus $\mathnormal{T}$ satisfies the right general maximal rank property. Any subspace $\mathnormal{V'}$ of $\mathnormal{V}$ of 
dimension 2 must have nonzero intersection with the span of 
$\{v_1,v_2\}$. Thus $\mathnormal{V'\otimes W}$ must have nonzero 
intersection with the span of $\{v_1\otimes w_1,v_2\otimes w_1\}$. 
This means that the restriction of $\mathnormal{T}$ to 
$\mathnormal{V'\otimes W}$ cannot have maximal rank.
\end{example}


\begin{thebibliography}{99}

\bibitem{ars}
M. Auslander, I. Reiten and S. O. Smal\o,
{\em Representation Theory of Artin Algebras},
Cambridge Studies in Advanced Mathematics, vol.~36, 
Cambridge University Press, New York, 1994.

\bibitem{G1}
J. A. Green,
\emph{Locally 
Finite Representations},
J. Alg. \textbf{41}~(1976), 
137-171.

\bibitem{gh}
P. Griffiths and J. Harris, {\em Principles of algebraic geometry}, Pure and Applied Mathematics, Wiley-Interscience [John Wiley \& Sons], New York, 1978.

\bibitem{HL1}
M. Hochster, D. Laksov,
\emph{The Linear 
Syzygies of Generic Forms},
Comm. Alg. \textbf{15}~(1987), no. 1-2, 
227-239.
\bibitem{K1}
M. Kleiner,
\emph{The Graded Preprojective Algebra of 
a Quiver},
Bull. London Math. Soc. \textbf{36}~(2004), 13-22.
\
\end{thebibliography}
\end{document}